\magnification=\magstep1
\hsize=16truecm
 
\input amstex
\TagsOnRight
\parindent=20pt
\parskip=2.5pt plus 1.2pt

\define\({\left(}
\define\){\right)}
\define\[{\left[}
\define\]{\right]}
\define\e{\varepsilon}

\define\supp {\sup\limits}

\define\summ{\sum\limits}

\define\bigcupp{\bigcup\limits}

\centerline{\bf On the tail behaviour of the distribution function
 of the maximum} 
\centerline{\bf for the partial sums of a class of i.i.d. random variables.}

\smallskip
\centerline{\it P\'eter Major}
\centerline{\it Alfr\'ed R\'enyi Mathematical Institute of the
Hungarian Academy of Science}
\centerline {e-mail address: major.peter$\@$renyi.mta.hu}

\medskip\noindent
{\bf Summary.} {\it We take an $L_1$-dense class of functions 
$\Cal F$ on a measurable space $(X,\Cal X)$ and a sequence of 
i.i.d. $X$-valued random variables $\xi_1,\dots,\xi_n$, and
give a good estimate on the tail behaviour of 
$\sup\limits_{f\in\Cal F}\sum\limits_{j=1}^nf(\xi_j)$ if the
conditions $\sup\limits_{x\in X}|f(x)|\le1$, $Ef(\xi_1)=0$
and $Ef(\xi_1)^2<\sigma^2$ with some $0\le\sigma\le1$ hold for 
all $f\in\Cal F$. Roughly speaking this estimate states that 
under some natural conditions the above considered supremum 
is not  much larger than the worst element taking part in it. 
The proof heavily depends on the main result of paper~[3]. 
Here we have to deal with such a problem where the classical 
methods worked out to investigate the behaviour of Gaussian 
or almost Gaussian random variables do not work.}

\beginsection 1. Introduction.

The main result of this paper is an estimate about the 
tail-distribution of the supremum of partial sums of
i.i.d. random variables presented in Theorem~1 together 
with an extension of it that provides an estimate 
for this tail-distribution in some cases not covered 
in Theorem~1. At first glance these results may look
rather complicated, but as I try to explain in Section~2
they yield sharp estimates under natural conditions. 
They express such a fact that under some natural 
conditions we can get an almost as good bound for the 
supremum of an appropriately defined class of partial 
sums as for one term taking part in this supremum. 
Before presenting these results I recall the definition
 of $L_1$-dense classes of functions, a notion that 
appears in the formulation of Theorem~1.

\medskip\noindent
{\bf Definition of $L_1$-dense classes of functions.}
{\it Let a measurable space $(X,{\Cal X})$ be given 
together with a class of ${\Cal X}$ measurable, real 
valued functions $\Cal F$ on this space. The class of 
functions ${\Cal F}$ is called an $L_1$-dense class of 
functions with parameter~$D$ and exponent~$L$ if for 
all numbers $0<\varepsilon\le1$ and probability 
measures $\nu$ on the space $(X,{\Cal X})$ there 
exists a finite $\varepsilon$-dense subset
${\Cal F}_{\varepsilon,\nu}=\{f_1,\dots,f_m\}\subset {\Cal F}$
in the space $L_1(X,{\Cal X},\nu)$ with 
$m\le D\varepsilon^{-L}$ elements, i.e. there exists 
such a set ${\Cal F}_{\varepsilon,\nu}\subset {\Cal F}$ 
with $m\le D\varepsilon^{-L}$ elements for which
$\inf\limits_{f_j\in{\Cal F}_{\varepsilon,\nu}}\int |f-f_j|\,d\nu
<\varepsilon$ for all functions $f\in {\Cal F}$.} 

\medskip
Theorem 1 yields the following estimate.

\medskip\noindent
{\bf Theorem 1.} {\it Let a sequence of independent, 
identically distributed random variables $\xi_1,\dots,\xi_n$, 
$n\ge2$, with values in a measurable space $(X,\Cal X)$ and 
with some distribution~$\mu$ be given together with a countable 
$L_1$-dense class of functions $\Cal F$ with parameter $D\ge1$ 
and exponent $L\ge1$ on the space $(X,\Cal X)$ such that 
$\supp_{x\in X}|f(x)|\le 1$, $\int f(x)\mu(\,dx)=0$, and 
$\int f^2(x)\mu(\,dx)\le\sigma^2$ with some number 
$0\le\sigma^2\le1$ for all $f\in\Cal F$. Define the normalized 
random sums  $S_n(f)=\frac1{\sqrt n}\summ_{j=1}^n f(\xi_j)$ 
for all $f\in\Cal F$. There are some universal constants 
$C_j>0$,  $1\le j\le 5$, (such that also the inequality $C_2<1$ 
holds), for which the inequality
$$
P\(\sup_{f\in\Cal F}|S_n(f)|\ge v\)\le 
C_1 e^{-C_2\sqrt n v\log(v/\sqrt n\sigma^2)}
\quad \text{for all } v\ge u(\sigma) \tag1.1
$$
holds if one of the following conditions is satisfied.

\medskip
\item{(a)} $\sigma^2\le\frac1{n^{200}}$, and  
$u(\sigma)=\frac{C_3}{\sqrt n}(L+\frac{\log D}{\log n})$,
\item{(b)} $\frac1{n^{200}}<\sigma^2\le\frac {\log n}{8n}$, 
and $u(\sigma)=\frac{C_4}{\sqrt n}
\(L\frac{\log n}{\log(\frac{\log n}{n\sigma^2})}+\log D\)$,
\item{(c)} $\frac{\log n}{8n}<\sigma^2\le1$, and
$u(\sigma)=\frac{C_5}{\sqrt n}(n\sigma^2+L\log n+\log D)$.} 

\medskip
I complete the result of Theorem~1 with an extension which 
is actually a repetition of Theorem~4.1 in~[2]. It yields an 
estimate for $P\(\sup\limits_{f\in\Cal F}|S_n(f)|\ge v\)$
in cases not covered in Theorem~1.

\medskip\noindent
{\bf Extension of Theorem 1.} {\it Let us consider, similarly
to Theorem~1, a sequence of independent, identically 
distributed random variables $\xi_1,\dots,\xi_n$, $n\ge2$, 
with values in a measurable space $(X,\Cal X)$ with some 
distribution $\mu$ together with a countable $L_1$-dense 
class of functions $\Cal F$ with parameter $D\ge1$ and 
exponent $L\ge1$ on the space $(X,\Cal X)$ such that 
$\supp_{x\in X}|f(x)|\le 1$, $\int f(x)\mu(\,dx)=0$, and 
$\int f^2(x)\mu(\,dx)\le\sigma^2$ with some number 
$0\le\sigma^2\le1$ for all $f\in\Cal F$. The supremum of the 
normalized partial sums $S_n(f)$, $f\in\Cal F$, introduced 
in Theorem~1 satisfies the inequality 
$$
P\left(\sup_{f\in{\Cal F}}|S_n(f)|\ge v\right)
\le C\exp\left\{-\alpha\frac {v^2}{\sigma^2}\right\} \tag1.2
$$
with appropriate (universal) constants $\alpha>0$, $C>0$ 
and $C_6>0$ if $\frac{\log n}{8n}<\sigma^2\le1$,
$\bar u(\sigma)\le v\le\sqrt n\sigma^2$, where 
$\bar u(\sigma)$ is defined as $\bar u(\sigma)
=C_6\sigma(L^{3/4}\log^{1/2}\frac2\sigma +(\log D)^{3/4})$.}

\medskip
The value $\frac{\log n}{8n}$ determining the boundary between 
cases~(b) and~(c) in Theorem~1 could by replaced by 
$\alpha\frac{\log n}n$ with any number $0<\alpha<1$. To see this
one has to check that the formula defining $u(\sigma)$ in
cases~(b) and~(c) give a value of the same order if 
$\sigma^2\sim \alpha\frac{\log n}n$ with $0<\alpha<1$. I chose
the parameter $\alpha=\frac18$ because some calculations were
simpler with such a choice. Let me remark that a similar 
statement holds for the value of boundary $n^{-200}$ between
cases~(a) and~(b). This could have been replaced by 
$n^{-\beta}$ with any $\beta>1$.

\medskip
In Section~2 I try to explain why the above results are natural,
in Section~3 I present their proof, and in Section~4 I make some
additional remarks. I finish this section with a short comparison 
of the results of this paper with some similar results of
Talagrand in~[6].

In both works the magnitude of the supremum of partial sums of
i.i.d. random variables are studied, and behind the results 
there is their implicit comparison with analogous estimates 
about the supremum of Gaussian random variables.

The analogous problems about the supremum of Gaussian random 
variables can be well investigated by the so-called chaining 
argument, whose best, sharpest version is worked out in~[6]. 
The estimates about the supremum of Gaussian random variables 
can be simply generalized for the supremum of other classes 
random variables if the tail distributions of the differences 
of the elements from these classes of random variables satisfy 
an estimate similar to the corresponding estimate in the 
Gaussian case. If we consider partial sums of independent 
random variables, then the tail distributions of these partial 
sums satisfy only a weaker estimate. Hence some additional 
conditions have to be imposed in order to get good results. 
Both here and in~[6] good estimates are given for the tail 
distribution of partial sums of i.i.d. random variables 
under some additional conditions. But these additional
conditions are different in the two works, and in my opinion 
the difference between them is not such a technical detail 
as it may seem at first sight. 

Talagrand extends the chaining argument to other models by
exploiting that under some additional conditions a better 
(Gaussian) estimate can be given for the tail distribution 
of sums i.i.d. random variables. His proof can be considered 
as the extension of a Gaussian argument to a more general 
class of models with `almost Gaussian behaviour'. The 
additional condition of this paper about the existence of 
an $L_1$-dense class of functions has a different character. 
It is useful to guarantee that the influence of some 
unpleasant `non-Gaussian effects' in the model we are 
working with is small. The proof of this fact demands an
argument different from the usual methods applied in the
Gaussian case. I do not write down the details about the 
difference of the two methods, because I did it in 
Chapter~18 of~[2] at pp.~235--237. Let me remark that 
the results obtained with their help cannot substitute 
each other. There are problems where the first one is 
useful and there are problems where the second one. 

\beginsection 2. Discussion on the conditions of these results.

Our goal was to give sharp estimate for the supremum of 
a class of normalized partial sums $S_n(f)$ defined in 
Theorem~1 if the functions $f$ are elements of an 
$L_1$-dense class of functions $\Cal F$ that satisfies 
the conditions of Theorem~1. We have to explain why 
formulas~(1.1) and~(1.2) provide the right estimate 
in this problem, and why we had to impose the conditions 
$v\ge u(\sigma)$ and $\bar u(\sigma)\le v\le\sqrt n\sigma^2$
in them. We prove such estimates which depend on some 
universal multiplying constants whose optimal choice we 
do not investigate. Besides, we try to give a good value 
for the functions $u(\sigma)$ and $\bar u(\sigma)$ only 
in the case when the parameter $D$ and exponent $L$ of 
the $L_1$-dense class of functions $\Cal F$ are bounded 
by a fixed number not depending on the parameter~$\sigma^2$. 
If the parameter~$D$ or exponent~$L$ is very large, then a 
different function~$u(\sigma)$ could be chosen that 
provides a sharper result.

If we disregard the value of the universal constants appearing
in our estimates then we can say that the estimate~(1.1) for 
the tail distribution of the supremum we consider and
the estimate of Bennett's inequality for the tail distribution
of a single term in this supremum agree, at least in the case 
if we consider the estimate of Bennett's inequality at level
$v\ge 2\sqrt n\sigma^2$. (This follows e.g. from formula~(3.3)
in this paper. We recalled Bennett's inequality, and  
formula~(3.3) is a part of it.) On the other hand, we 
considered in Theorem~1 only such levels~$v$ where this 
condition is satisfied, since $u(\sigma)\ge2\sqrt n\sigma$ 
in all cases of Theorem~1. Moreover, there are examples that 
show that inequality~(3.3) is sharp, we cannot get a better 
estimate without some additional restrictions. (See Example~3.3 
in~[2]). The estimate~(1.2) in the extension of Theorem~1 in 
the case $\bar u(\sigma)\le v\le\sqrt n\sigma^2$ is also sharp 
(we disregard again the value of the universal constants in 
this formula), since the tail-distribution of a normalized 
partial sum cannot have a better bound, than the Gaussian 
estimate given in~(1.2). Formally there is a gap between the 
results of Theorem~1 and its extension, because we did not 
consider the case $\sqrt n\sigma^2\le v\le u(\sigma)$. But 
this gap can be simply filled in the case when the numbers~$D$ 
and~$L$ are bounded by constants not depending on~$\sigma^2$, 
and we do not try to find optimal universal constants in our 
estimates. Indeed, in this case we have 
$u(\sigma)\le\frac{\bar C}{\sqrt n}n\sigma^2$, and
$$
P\(\sup_{f\in\Cal F}|S_n(f)|\ge v\)\le 
P\(\sup_{f\in\Cal F}|S_n(f)|\ge \sqrt n\sigma^2\)\le 
Ce^{-\alpha n\sigma^2}\le Ce^{-\bar\alpha v^2/\sigma^2},
$$ 
i.e. relation~(1.2) holds (with a possible different parameter
$\bar\alpha>0$) for all $\bar u(\sigma)<v\le u(\sigma)$. This
estimate is sharp again.

We also have to understand why we could give a good estimate 
for the supremum of normalized partial sums only under the
conditions $v\ge u(\sigma)$ in cases~(a) and~(b) and 
$v\ge\bar u(\sigma)$ in case~(c). I shall present an example
that satisfies the conditions of Theorem~1, and in which 
there is no useful estimate in formulas~(1.1) and~(1.2) for 
$v<u(\sigma)$ in cases~(a) and~(b) or $v<\bar u(\sigma)$ in
cases~(c). (More precisely, we allow a different multiplying 
factor $C_j$ as in the definition of $u(\sigma)$ and 
$\bar u(\sigma)$ when we consider this model.) This implies 
in particular that the conditions $v\ge u(\sigma)$ and 
$v\ge\bar u(\sigma)$ cannot be dropped in Theorem~1 and 
in its extension. 

At this point it may be useful to recall the concentration 
inequality for the supremum of partial sums of independent 
random variables. (See e.g.~[5]). It states that there is a 
concentration point of the supremum of partial sums of 
independent random variable such that this supremum is 
strongly concentrated in a small neighbourhood of this 
concentration point. I do not formulate this result in a 
more precise and detailed form, because we need it here 
only for the sake of some orientation. The problem with 
the application of this result is that it determines the 
concentration point only in an implicit way as the expected 
value of the supremum we are investigating, and we cannot 
calculate it explicitly in the general case. On the other 
hand, the concentration inequality implies that we can get
a good, non-trivial estimate for the tail distribution of 
partial sums of independent random variables only at levels 
higher than the concentration point of the partial sums. 
(We call such estimates trivial which only say that a 
probability is not greater than~1.) So the numbers 
$u(\sigma)$ and $\bar u(\sigma)$ in Theorem~1 and in its 
extension are actually upper bounds for the concentration 
point of the supremum, and we shall present a model 
satisfying the conditions of Theorem~1, where the values 
$u(\sigma)$ and $\bar u(\sigma)$ determine the concentration 
point of the supremum up to a multiplicative factor.

We shall consider the following model. Take independent, 
uniformly distributed random variables $\xi_1,\dots,\xi_n$ 
on the unit interval~$[0,1]$, fix a number $0\le\sigma^2\le 1$, 
and define a class of functions $\Cal F_\sigma$ and 
$\bar {\Cal F_\sigma}$ with functions defined on the unit interval 
$[0,1]$ in the following way. $\Cal F_\sigma=\{f_1,\dots,f_k\}$, 
and $\bar{\Cal F}=\{\bar f_1,\dots,\bar f_k\}$ with 
$k=k(\sigma)=[\frac1{\sigma^2}]$, where $[\cdot]$ denotes 
integer part, and $\bar f_j(x)=\bar f_j(x|\sigma)=1$  if 
$x\in[(j-1)\sigma^2,j\sigma^2)$, $\bar f_j(x)=\bar f_j(x|\sigma)=0$ 
if  $x\notin[(j-1)\sigma^2,j\sigma^2)$, $1\le j\le k$, and 
$f_j(x)=f_j(x|\sigma)=\bar f_j(x)-\sigma^2$, $1\le j\le n$. 
It can be seen that $\Cal F_\sigma$ satisfies the conditions 
of Theorem~1 with parameter~$\sigma^2$. In particular, it 
is an $L_1$-dense class with such a parameter~$D$ and 
exponent~$L$ that can be bounded by numbers not depending 
on~$\sigma^2$. This can be seen directly, but it is also a 
consequence of some classical results by which the indicator 
functions of a Vapnik--\v{C}ervonenkis class of sets constitute 
an $L_1$~dense class of functions. (See e.g. Theorem~5.2 in~[2]).

I shall show that in this example a number $\bar C>0$ can be 
chosen in such a way that for all $\delta>0$ there is an index 
$n_0(\delta)$ such that for all sample sizes 
$n\ge n_0(\delta)$ and numbers $0\le\sigma\le1$ the inequality 
$$
P\left(\sup_{f\in{\Cal F_\sigma}}|S_n(f)|
\ge\hat u(\sigma)\right)\ge1-\delta, \tag2.1
$$
holds with $\hat u(\sigma)=\frac{\bar C}{\sqrt n}$ in case~(a),
$\hat u(\sigma)=\frac{\bar C}{\sqrt n}\frac{\log n}
{\log(\frac{\log n}{n\sigma^2})}$ in case~(b), and
$\hat u(\sigma)=\bar C\sigma\log^{1/2}\frac2\sigma$ in case~(c).
This result may explain why we had to impose the conditions 
$v>u(\sigma)$ and $v>\bar u(\sigma)$ in Theorem~1 and in its 
extension. (We are interested only in such cases when the
estimate of Theorem~1 or its extension provide an upper bound
strictly less than~1, (i.e. smaller than a number $\alpha<1$ 
for all parameters $\sigma^2$ and~$n$), and this is the case 
if the constants $C_j$, $j=3,4,5,6$  are chosen sufficiently 
large in these results.) 

To prove relation~(2.1) introduce the following notation.
Define the empirical distribution function $F_n(x)$ of the
random variables $\xi_1,\dots,\xi_n$, i.e. put
$$
F_n(x)=\frac1n\{\text{the number of indices }j,\; 1\le j\le n,
\text{ such that } \xi_j< x\}
$$ 
for all $0<x\le 1$, and take its normalization 
$G_n(x)=\sqrt n(F_n(x)-x)$, $0<x\le 1$. Observe that
$$
\left\{\sup_{f\in{\Cal F_\sigma}}|S_n(f)|\ge\hat u(\sigma)\right\}
=\left\{\max_{1\le j\le k(\sigma)}|G_n(j\sigma)-G_n((j-1)\sigma)|
\ge\hat u(\sigma)\right\}. \tag2.2
$$
By a classical results of probability theory, the normalized
empirical distribution functions weakly converge to the Brownian
bridge as $n\to\infty$. In our next considerations it will be 
also interesting that the modulus of continuity of a Brownian 
bridge, (which actually agrees with the modulus of continuity 
of a Wiener process) can be also calculated. (see~e.g.~[4]). 
By a similar, but simpler calculation we can estimate the 
probability of the event we get by replacing the normalized 
empirical distribution function $G_n(\cdot)$ by a Brownian bridge 
in the right-hand side expression of~(2.2). This is actually done 
with the choice $u(\sigma)=\bar C\sigma\log^{1/2}\frac2\sigma$ 
in the fourth chapter of~[2] (page~27), and it is shown 
that this probability is almost one for large parameters~$n$
for all $\sigma>0$ if the coefficient $\bar C$ of $u(\sigma)$ 
is chosen sufficiently small. (Actually we have to choose
$\bar C<\sqrt 2$.) Let us call this estimate the Gaussian version
of formula~(2.1). At a heuristic level this result together 
with formula~(2.2) and the weak convergence of the normalized 
empirical processes $G_n(\cdot)$ to a Brownian bridge suggests
that formula~(2.1) should hold with 
$u(\sigma)=\bar C\sigma\log^{1/2}\frac2\sigma$  and a small 
coefficient $\bar C>0$. 

This heuristic argument is nevertheless misleading, since the
weak convergence of the empirical processes $G_n(\cdot)$ to
the Brownian bridge does not make possible to carry out a 
limiting procedure that leads to formula~(2.1). On the other 
hand, a stronger version of the weak convergence of the 
normalized empirical processes (see~[1]) yields a useful result 
in this direction. This result states a normalized empirical 
process~$G_n(x)$ and a Brownian bridge $B(x)$, $0\le x\le1$, 
can be constructed in such a way that 
$\sup\limits_{0\le x\le 1}|B(x)-G_n(x)|\le K\frac{\log n}{\sqrt n}$ 
for all $n\ge2$ and sufficiently large $K>0$ with probability 
almost~1. This result together with the Gaussian version of 
formula~(2.1) imply the validity of formula~(2.1) if 
$\sigma^2\ge B\frac{\log n}{2n}$ with a sufficiently 
large $B>0$. Indeed, in this case 
$\hat u(\sigma)\ge2K\frac{\log n}{\sqrt n}$, hence the
Gaussian version of formula of~(2.1) together with the 
result of~[1] imply that
$$P
\left(\max_{1\le j\le k(\sigma)}|G_n(j\sigma)-G_n((j-1)\sigma)|
\ge\frac{\hat u(\sigma)}2\right)\ge1-\delta
$$ 
if 
$\sigma^2\ge B\frac{\log n}n$, and $n\ge n_0(\delta)$, i.e. 
inequality~(2.1) holds in this case if we replace $\bar C$ by
$\frac{\bar C}2$ in the definition of~$\hat u(\sigma)$. Moreover,
this relation holds for all $\sigma^2\ge\frac{\log n}{8n}$, i.e.
in the case~(c) generally if we choose 
$\hat u(\sigma)=\bar C\sigma\log^{1/2}\frac2\sigma$ with a 
sufficiently small $\hat C>0$. To see this it is enough to
observe that if
$\max\limits_{1\le j\le k(\sigma)}|G_n(j\sigma)-G_n((j-1)\sigma)|
\le\hat u(\sigma)$, then for any positive integers~$A$ we have
$\max\limits_{1\le j\le k(\sqrt A\sigma)}
|G_n(j(A\sigma))-G_n((j-1)(A\sigma))|\le A\hat u(\sigma)$, 
and that the corresponding result holds if 
$\sigma^2\ge B\frac{\log n}{8n}$.

In cases (a) and (b) the above Gaussian approximation argument
does not work. In case~(b) we shall prove formula~(2.1) 
by means of a Poissonian approximation method described below. 
It can be considered as a more detailed elaboration of the 
argument in Example~4.3 of~[2].

In this argument first we consider the following problem. Take 
a Poisson process $Z_n(t)$, $0\le t\le 1$, with parameter $n$, 
(i.e. let $EZ_n(t)=nt$ for all $0\le t\le 1$) in the interval 
$[0,1]$. Fix some number $0\le\sigma^2\le\frac17\frac{\log n}n$, 
and define with its help the number
$\hat u(\sigma)=\hat u(\sigma,n)=\frac3{4\sqrt n}\frac{\log n}
{\log(\frac{\log n}{n\sigma^2})}$ and the random variables 
$\bar V_j=\bar V_j^{(n)}(\sigma)=Z_n(j\sigma^2)-Z_n((j-1)\sigma^2)$ 
for $1\le j\le k$ with $k=[\frac1{\sigma^2}]$. (Here we defined 
$\hat u(\sigma)$ similarly to quantity introduced with the same
notation at the formulation of inequality~(2.1) in the case~(b). 
We only made small modifications. Namely we considered $\sigma^2$ 
in the interval $[0,\frac17\frac{\log n}n]$ instead of the interval 
$[\frac1{n^{200}},\frac{\log n}{8n}]$, and we fixed the value
$\bar C=\frac34$ in the definition of $\hat u(\sigma)$. We want to 
show that for all $\delta>0$ there is some threshold index 
$n_0(\delta)$ such that the inequality
$$
P\left(\max_{1\le j\le k(\sigma)}\bar V_j^{(n)}(\sigma)\ge
\sqrt n\hat u(\sigma,n)\right)\ge 1-\delta
\quad \text{if } n\ge n_0(\delta) \tag2.3
$$
holds for all $0\le\sigma^2\le\frac17\frac{\log n}n$.

To prove this inequality let us first observe that
$$
\align
P\left(\max_{1\le j\le k(\sigma)}\bar V_j^{(n)}(\sigma)\ge
\sqrt n\hat u(\sigma,n)\right)
&\ge P(\bar V_j^{(n)}(\sigma)=
\sqrt n\hat u(\sigma,n)\text{ for some }1\le j\le k) \\
&=1-P(\bar V_1^{(n)}(\sigma)\neq
\sqrt n\hat u(\sigma,n))^k,
\endalign
$$
and
$$
\align
&P(\bar V_1^{(n)}(\sigma)\neq\sqrt n\hat u(\sigma,n))
=1-P(\bar V_1^{(n)}(\sigma)=\sqrt n\hat u(\sigma,n)) \\
&\qquad =1-\frac{(n\sigma^2)^{\sqrt n\hat u(\sigma,n)}}
{(\sqrt n\hat u(\sigma,n))!}e^{-n\sigma^2}
\le1-\left(\frac{n\sigma^2}
{\sqrt n\hat u(\sigma,n)}\right)^{\sqrt n\hat u(n,\sigma)}
e^{-n\sigma^2}.
\endalign
$$
Since we have $k=\frac 1{\sigma^2}$ we can bound the left-hand side
of~(2.3) from below as
$$
P\left(\max_{1\le j\le k(\sigma)}\bar V_j^{(n)}(\sigma)
\sqrt n\hat u(\sigma,n)\right)
\ge1-\left[1-\left(\frac{n\sigma^2}{\sqrt n\hat u(\sigma,n)}\right)
^{\sqrt n\hat u(n,\sigma)}e^{-n\sigma^2}\right]^{1/\sigma^2}
\ge1-e^{-T}
$$
with $T=\frac1{\sigma^2}
\left(\frac{n\sigma^2}{\sqrt n\hat u(\sigma,n)}\right)
^{\sqrt n\hat u(n,\sigma)}e^{-n\sigma^2}$, hence to prove~(2.3) it 
is enough to show that
$$
\left(\frac{n\sigma^2}{\sqrt n\hat u(\sigma,n)}\right)
^{\sqrt n\hat u(n,\sigma)}\ge\sigma^2e^{n\sigma^2}\log\frac1\delta 
\quad \text{if }n\ge n_0(\delta).
\tag2.4
$$

The right-hand side of~(2.4) can be bounded from above as
$$
\sigma^2e^{n\sigma^2}\log\frac1\delta 
=\frac{\log\frac1\delta}n (n\sigma^2)e^{n\sigma^2}
\le\frac{\log\frac1\delta}n \left(\frac17\log n\right)
e^{(\log n)/7}\le n^{-5/6}
$$
if $n\ge n_0(\delta)$, since $n\sigma^2\le\frac17\log n$. Hence
we prove~(2.4) if we show that
$$
\frac{\sqrt n\hat u(n,\sigma)}{n\sigma^2}
\log\left(\frac{\sqrt n\hat u(\sigma,n)}{n\sigma^2}\right)
\le\frac56\frac{\log n}{n\sigma^2}.
$$
By applying the definition of $\hat u(n,\sigma)$ and introducing
the quantity $z=\frac34\frac{\log n}{n\sigma^2}$ we can rewrite the
last inequality as 
$\frac z{\log (\frac{4z}3)}\log(\frac z{\log (\frac{4z}3)})\le\frac{10}9z$,
or since $z\ge\frac{21}4$ in the case we are investigating it can be 
rewritten as 
$\frac19\log\frac{4z}3\ge -\log\log\frac{4z}3-\log\frac43$ if 
$z\ge\frac{21}4$, and this relation clearly holds. Thus we proved~(2.3).

We shall prove relation~(2.1) in the case~(b) by means of 
formula~(2.3) for a Poisson process with parameter $\frac{99}{100}n$ 
instead of~$n$ and a simple coupling argument between an empirical 
process and a Poisson process. Namely, we make the following 
coupling. Let us consider a sequence of independent random variables 
$\xi_1,\xi_2\dots$ with uniform distribution on the unit interval
$[0,1]$ together with a Poissonian random variable $\eta=\eta_n$ 
with parameter $\frac{99}{100}n$ independent of the random
variables $\xi_j$, $j=1,2,\dots$, and take the first $\eta_n$ terms
of the random variables $\xi_j$, i.e. the sequence
$\xi_1,\xi_2,\dots,\xi_{\eta_n}$ with the random stopping index 
$\eta_n$. In such a way we constructed a Poisson process with 
parameter $\frac{99}{100}n$, which is smaller than the (non-normalized)
empirical distribution of the sequence $\xi_1,\dots,\xi_n$ in the
following sense. For large parameter $n$ with probability almost~1
all intervals $[a,b]\subset[0,1]$ contain more points from the 
sequence $\xi_1,\dots,\xi_n$ than from the above constructed
Poisson process. This is a simple consequence of the fact that
$P(\eta_n>n)\to0$ as $n\to\infty$.

The above coupling construction and formula~(2.3) (with a Poisson
process with parameter $\frac{99}{100}$) imply that
$$
P\left(\sup_{\bar f\in\bar{\Cal F}_\sigma}\sqrt nS_n(\bar f)\ge
\sqrt{\frac{99}{100} n}
\hat u(\sigma,\left(\frac{99}{100}n\right)\right)\ge 1-\delta
\quad \text{if } n\ge n_0(\delta) 
$$
with the class of functions $\bar{\Cal F}_\sigma$ introduced before
the formulation~(2.1) and the function $\hat u(\sigma,n)$ defined in
the discussion of case~(b). To complete the proof of~(2.1) in the 
case~(b) it is enough to check that the above relation remains valid
if the class of functions $\bar{\Cal F}_\sigma$ is replaced by the
class of functions $\Cal F_\sigma$ and the term
$\sqrt{\frac{99}{100} n}\hat u(\sigma,\frac{99}{100}n)$ is replaced by
$\hat u(\sigma,n)=\frac{\bar C}{\sqrt n}\frac{\log n}
{\log(\frac{\log n}{n\sigma^2})}$ with some appropriate $\bar C>0$.
Since the functions $f\in\Cal F$ are of the form 
$f(x)=\bar f(x)-\sigma^2$ with some $\bar f\in\Cal F$, this has the
consequence $\sqrt nS_n(f)=\sqrt nS_n(\bar f)-n\sigma^2$, and to 
prove the desired relation it is enough to check that 
$$
\sqrt{\frac{99}{100} }\frac34\frac{\log n}
{\log(\frac{\log n}{\frac{99}{100}n\sigma^2})}-n\sigma^2\ge 
\sqrt{\frac{99}{100} }\frac34\frac{\log n}
{\log(\frac{\log n}{n\sigma^2})}-n\sigma^2\ge 
\bar C \frac{\log n}{\log(\frac{\log n}{n\sigma^2})}
$$
with some appropriate $\bar C>0$ if $8n\sigma^2\le\log n$. The 
first inequality clearly holds, and the second inequality is 
equivalent to the relation
$$
\sqrt{\frac{99}{100} }\frac34\frac{\frac{\log n}{n\sigma^2}}
{\log(\frac{\log n}{n\sigma^2})}\ge \alpha
$$
with some $\alpha>1$. But this relation clearly holds if 
$8n\sigma^2\le\log n$. Thus we have proved~(2.1) also in case~(b).

In the case~(a) the proof of~(2.1) is very simple. It is enough
to observe that the sample points $\xi_j$ fall into one of the 
intervals $[(j-1)\sigma^2,j\sigma^2)$, $1\le j\le k$, (we disregard
the event that they fall into the last interval $[k\sigma^2,1)$ which
has negligible small probability), hence
$$
P\left(\sup_{\bar f\in\bar{\Cal F}_\sigma}\sqrt nS_n(\bar f)=1\right)
\ge 1-\delta \quad \text{if } n\ge n_0(\delta), 
$$
and since $\sigma^2$ is very small for large $n$ relation~(2.1)
holds in case~(a) with $\bar C=1-\e$ for any $\e>0$.

\medskip
At the end of this section let me remark that in the above 
example actually we have given a lower bound on the modulus 
of continuity of a normalized empirical process. I formulate 
a problem below where the proof of a stronger form of this 
result is suggested.

\medskip\noindent
{\bf Problem.} {\it Let $\xi_1,\xi_2,\dots$ be a sequence of 
independent random variables, uniformly distributed in the unit
interval $[0,1]$, and define with its help the empirical 
distribution functions 
$$
F_n(x)=\frac1n\text{ times the number of indices $j$, $1\le j\le n$,
such that } \xi_j<x
$$ 
together with their normalizations $G_n(x)=\sqrt n(F_n(x)-x)$, 
$0\le x\le 1$, for all indices $n=1,2,\dots$. Find such a 
function $v(n,\sigma^2)$, $n=1,2,\dots$, $0\le\sigma^2\le1$, 
for which
$$
\lim_{n\to\infty}\sup_{\{(s,t)\colon\; 0\le s,t\le1,\; |t-s|\le\sigma_n^2\}}
\frac{|G_n(t)-G_n(s)|}{v(n,\sigma^2_n)}=1\quad \text{with probability }1
$$
if $\sigma^2_n\to0$ as $n\to\infty$.}

\beginsection 3. Proof of Theorem 1 and its extension.

{\it Proof of Theorem 1.}\/ In the case (a) inequality (1.1) is a simple
consequence of Theorem~1 in~[3]. We can apply this result (by writing
$\sigma^2$ instead of $\rho$ in its formulation), since 
$\int f^2(x)\mu(\,dx)\le \int|f(x)|\mu(\,dx)$ if 
$\supp_{x\in X}|f(x)|\le1$, hence under the conditions of Theorem~1 
the inequality $\int|f(x)|\mu(\,dx)\le\rho$ holds for all 
$f\in{\Cal F}$ with $\rho=\sigma^2$. Hence
$$
P\left(\sup_{f\in{\Cal F}}|S_n(f)|\ge v\right)
\le De^{-\frac1{50} \sqrt nv\log(\sigma^{-2})}
\quad\text{if } v\ge\frac{\bar C}{\sqrt n}L \text{ and }
\sigma^2\le\frac1{n^{200}} \tag3.1
$$
with an appropriate $\bar C>0$.

I claim that we can drop the coefficient $D$ at the right-hand 
side of~(3.1) if we replace the coefficient $\frac1{50}$ by
$\frac1{100}$ in the exponent, we choose such a constant $\bar C$
in~(3.1) for which $\bar C\ge\frac12$, and impose condition~(a) in 
the form $v\ge\frac{\bar C}{\sqrt n}(L+\frac{\log D}{\log n})$. 
To show this it is enough to check that 
$D\le e^{\frac1{100}\sqrt nv\log(\sigma^{-2})}$
in this case. This relation holds, since 
$\frac{\log D}{\log n}\le 2\sqrt nv$, and 
$\log(\sigma^{-2})\ge 200\log n$, thus
$D=\exp\{\frac1{200}(\frac{\log D}{\log n})(200\log n)\}
\le\exp\{\frac1{100}\sqrt nv\log(\sigma^{-2})\}$, as I claimed.

Next I show that formula~(3.1) or its previous modification
remains valid if we replace $\log(\sigma^{-2})$ by 
$\log(\frac v{\sqrt n\sigma^2})$ in the exponent of its 
right-hand side. In the proof of this statement we can 
restrict our attention to the case $v\le\sqrt n$, since 
otherwise the probability at the left-hand side of~(3.1) 
equals zero. In this case the inequality 
$\sigma^{-2}\ge\frac v{\sqrt n\sigma^2}$ holds, and this 
allows the above replacement. The above modifications of 
formula~(3.1) imply inequality~(1.1) in case~(a).

\medskip\noindent
{\it Remark.}\/ If we are not interested in the value of the 
(universal) constants in~(1.1), then in the case~(a) this 
inequality has the same strength if we replace the term 
$\log (v/\sqrt n\sigma^2)$ by $\log(\sigma^{-2})$ in it. To 
see this, observe that beside the inequality 
$\sigma^{-2}\ge\frac v{\sqrt n\sigma^2}$ (if $v\le\sqrt n$), 
the inequality $\frac v{\sqrt n\sigma^2}\ge\frac1{n\sigma^2}
\ge\sigma^{-2+1/100}$ also holds in case~(a) because 
of the inequalities $v\ge u(\sigma)\ge n^{-1/2}$ and 
$n^{-200}\ge\sigma^2$. The original form of~(1.1) has the 
advantage that it simultaneously holds in all cases~(a), (b) 
and~(c).

\medskip\noindent
{\it The proof of Theorem~1 in cases~(b) and~(c)}.\/
By applying the $L_1$-dense property of the class of functions $\Cal F$
with the choice $\e=n^{-200}$ and the measure $\mu$ we may find 
$m\le Dn^{200L}$ functions $f_j\in\Cal F$, $1\le j\le m$, such that
$\min\limits_{1\le j\le m}\int |f_j(x)-f(x)|\mu(\,dx)\le n^{-1/200}$
for all $f\in\Cal F$. This means that $\Cal F=\bigcupp_{j=1}^n\Cal D_j$
with 
$$
\Cal D_j=\left\{f\colon\; f\in\Cal F,\,\int|f_j(x)-f(x)|\mu(\,dx)\le n^{-200}
\right\},
$$
and as a consequence
$$
P\(\sup_{f\in\Cal F}|S_n(f)|\ge v\)\le \sum_{j=1}^m
P\(|S_n(f_j)|\ge\frac v2\)
+\sum_{j=1}^m P\(\sup_{f\in\Cal D_j}|S_n(f-f_j)|\ge\frac v2\) \tag3.2
$$
for all $v>0$. We shall estimate both terms at the right-hand side of~(3.2)
if $v\ge u(\sigma)$, the first one by means of Bennett's inequality,
more precisely by a consequence of this inequality formulated below,
and the second term by means of the already proved case~(a) of Theorem~1.
We shall apply the following version of Bennett's inequality, see~[2].

\medskip\noindent
{\bf Bennett's inequality.} {\it Let $X_1,\dots,X_n$ 
be independent and identically distributed random variables such
that, $P(|X_1|\le1)=1$, $EX_1=0$, and $EX_1^2\le\sigma^2$ with
some $0\le\sigma\le 1$. Put 
$S_n=\frac1{\sqrt n}\sum\limits_{j=1}^n X_j$. Then
$$
P(S_n>v)\le\exp\left\{-n\sigma^2\left[\left(1+\frac v{\sqrt n\sigma^2}\right)
\log\left(1+\frac v{\sqrt n\sigma^2}\right)
-\frac v{\sqrt n\sigma^2}\right]\right\}
\quad\text{for all } v>0. 
$$
As a consequence, for all $\varepsilon>0$ there exists some
$B=B(\varepsilon)>0$ such that
$$
P\left(S_n>v\right)\le\exp\left\{-(1-\varepsilon)\sqrt nv
\log \frac v{\sqrt n\sigma^2}
\right\}\quad\text{if } v>B\sqrt n\sigma^2, 
$$
and there exists some positive constant $K>0$ such that
$$
P\left(S_n>v\right)\le\exp\left\{-K\sqrt nv\log \frac v{\sqrt n\sigma^2}
\right\}\quad\text{if }v>2\sqrt n\sigma^2. \tag3.3
$$
}

\medskip
The above result is a special case of Theorem~3.2 in~[2], in the case
when we restrict our attention to sums of independent and identically
distributed random variables. It has a slightly different form, because
in the definition of $S_n$ we considered normalized sums (with a 
multiplication by $n^{-1/2}$). Here we need only the inequality
formulated in~(3.3) which helps to estimate the probabilities 
appearing in the first sum at the right-hand side of~(3.2). 
To apply~(3.3) in the estimation of these terms we have to show that  
$u(\sigma)>2\sqrt n\sigma^2$ in cases~(b) and~(c) if the constants 
$C_4$ and $C_5$ are chosen sufficiently large in Theorem~1.

In case~(b) it is enough to show that $u(\sigma)\ge 
C_4\frac{\log n}{\log(\frac{\log n}{n\sigma^2})}\ge 2n\sigma^2$,
and even
$C_4\frac{\log n}{\log(\frac{\log n}{n\sigma^2})}\ge 20n\sigma^2$,
or in an equivalent form $\frac{C_4}{20}\frac{\log n}{n\sigma^2}
\ge\log(\frac{\log n}{n\sigma^2})$. (Observe that
$\frac{\log n}{n\sigma^2}\ge8$, hence 
$\log(\frac{\log n}{n\sigma^2})>0$ in case~(b).)
 This statement holds, since $z=\frac{\log n}{n\sigma^2}\ge2$ 
in case~(b), and $\frac{C_4}{20}z\ge\log z$ if $z\ge8$, and $C_4$ 
is sufficiently large.

In case~(c), clearly $u(\sigma)\ge\frac{C_5}{\sqrt n}n\sigma^2
\ge20\sqrt n\sigma^2$ for sufficiently large constant $C_5$. 
These relations together with formula~(3.3) imply that in 
cases (b) and (c)
$$
P\(|S_n(f_j)|\ge\frac v2\)\le 
2\exp\left\{-K\sqrt nv\log \frac v{\sqrt n\sigma^2}\right\} 
\quad\text{if }v\ge u(\sigma) \tag3.4
$$
with an appropriate $K>0$ for all $1\le j\le m$. (In formula~(3.4)
we exploit that $\log(\frac{\frac v2}{\sqrt n\sigma^2})\ge
\frac12\log(\frac v{\sqrt n\sigma^2})$ since 
$\frac v{\sqrt n\sigma^2}\ge20$, and as a consequence 
$\log(\frac v{\sqrt n\sigma^2})\ge 2\log 2$.)

Let us define, with the help of the class of functions $\Cal D_j$
the class of functions 
$\Cal D_j'=\{h\colon\; h=\frac{f-f_j}2,\;f\in\Cal D_j\}$ for
all $1\le j\le m$. It is not difficult to see that 
$\supp_{x\in X}|h(x)|\le1$, 
$\int h^2(x)\mu(\,dx)\le\int |h(x)|\mu(\,dx)\le n^{-200}$ for all
$h\in\Cal D'_j$, and $\Cal D_j'$ is an $L_1$-dense class of functions
with parameter $D$ and exponent $L$, $1\le j\le m$. I claim that
$$
\aligned
P\(\sup_{f\in\Cal D_j}|S_n(f-f_j)|\ge\frac v2\) 
&=P\(\sup_{h\in\Cal D'_j}|S_n(h_j)|\ge\frac v4\) \\
&\le e^{-C_2\sqrt nv\log(vn^{195})} \quad \text{if }v\ge u(\sigma)
\endaligned\tag3.5
$$  
for all $1\le j\le n$ in both cases~(b) and~(c). We shall get 
this estimate by applying Theorem~1 in the already proved 
case~(a) with the choice of parameter $\sigma^2=n^{-200}$. To 
apply this result we have to check that 
$\frac{u(\sigma)}4\ge u(n^{-200})
=\frac{C_3}{\sqrt n}(L+\frac{\log D}{\log n})$
if the constants $C_4$ and $C_5$ are sufficiently large. These 
statements hold, since in case~(b) 
$\frac{\log n}{\log \frac{\log n}{n\sigma^2}}
\ge\frac{\log n}{\log (n^{199}\log n)}\ge\frac1{200}$, hence 
$u(\sigma)\ge \frac{C_4}{\sqrt n}(\frac{L\log n}{200}+\log D)
\ge4u(n^{-200})$ if $C_4$ is chosen sufficiently large,
and an analogous but simpler argument supplies this relation
in case~(c) if $C_5$ is chosen sufficiently large.

It is not difficult to see that the right-hand side both of~(3.4)
and~(3.5) can be bounded from above by 
$C_1e^{-\bar C_2 \sqrt nv\log(v/\sqrt n\sigma^2)}$
with some appropriate constants $C_1>0$ and $\bar C_2>0$. Hence 
relations (3.2), (3.4) and (3.5) together with the inequality 
$m\le Dn^{200L}$ imply that
$$
P\(\sup_{f\in\Cal F}|S_n(f)|\ge v\)\le C_1 Dn^{200L}
e^{-\bar C_2 \sqrt nv\log(v/\sqrt n\sigma^2)} \quad\text{if }v\ge u(\sigma) 
\tag3.6
$$
in both cases (b) and (c). Hence to complete the proof of Theorem~1
(with the choice $C_2=\frac{\bar C_2}2$) it is enough to show that
$$
e^{-\frac{\bar C_2}2 \sqrt nv(\sigma)v/\sqrt n\sigma^2)}\le 
e^{-\frac{\bar C_2}2 \sqrt nu(\sigma)\log(u(\sigma)/\sqrt n\sigma^2)} 
\le D^{-1}n^{-200L} \quad\text{if }v\ge u(\sigma) \tag3.7
$$
in cases~(b) and~(c) if the constants $C_4$ and $C_5$ are chosen
sufficiently large.

It is enough to prove the second inequality in formula~(3.7),
since its proof also implies that the expressions in the 
exponent of this formula have negative value, and they 
are decreasing functions for $v\ge u(\sigma)$. The second
inequality in~(3.7) clearly holds in case~(c), since 
$\frac{\bar C_2}2\sqrt nu(\sigma)\ge200L\log n+\log D$,
and $\frac{\log(u(\sigma)}{\sqrt n\sigma^2}\ge1$ in this case. 
In case~(b) relation~(3.7) can be reduced to the inequalities
$\bar C_2 \sqrt nu(\sigma)\log(\frac{u(\sigma)}{\sqrt n\sigma^2})
\ge800L\log n$, and
$\bar C_2 \sqrt nu(\sigma)\log(\frac{\sqrt n(\sigma)}{n\sigma^2})
\ge4\log D$. To prove the second inequality observe that
in case~(b) 
$$
\bar C_2\sqrt nu(\sigma)\ge C_4\bar C_2\log D\ge4\log D, 
\quad \text{and} \quad\log\(\frac{\sqrt nu(\sigma)}{n\sigma^2}\)\ge1.
$$
The second of these inequalities follows from the relation
$\frac{\sqrt nu(\sigma)}{n\sigma^2}\ge C_4
\frac{\frac{\log n}{n\sigma^2}}{\log(\frac{\log n}{n\sigma^2})}\ge3$,
which holds because of the relation $\frac{\log n}{n\sigma^2}\ge8$
in case~(b).

The remaining inequality can be rewritten as
$\bar C_2 \frac{\sqrt nu(\sigma)}{n\sigma^2}\log(\frac{\sqrt nu(\sigma)}
{n\sigma^2})\ge800L\frac{\log n}{n\sigma^2}$. To prove it observe 
that because of the definition of the function $u(\sigma)$ in 
case~(b) we can write
$\bar C_2\frac{\sqrt nu(\sigma)}{n\sigma^2}\ge1600L\frac{\log n}{n\sigma^2}
\frac1{\log(\frac{\log n}{n\sigma^2})}\ge1600L\frac{\log n}{n\sigma^2}$,
since $\log(\frac{\log n}{n\sigma^2})\ge\log8\ge1$. I also claim that
$\log(\frac{\sqrt nu(\sigma)}{n\sigma^2})
\ge\frac12(\frac{\log n}{n\sigma^2})$. By multiplying the last two
inequalities we get the desired inequality, and this completes the 
proof of Theorem~1.

To prove the above formulated inequality introduce the notation 
$z=\frac{\log n}{n\sigma^2}$. By exploiting the definition of 
$u(\sigma)$ in case~(b) we can write with the help of this 
notation that
$\log(\frac{\sqrt nu(\sigma)}{n\sigma^2})\ge \log z-\log\log z
\ge\frac12\log z=\frac12(\frac{\log n}{n\sigma^2})$. In the 
above argument we have exploited that in case~(b) $z\ge8$, 
hence $\log z-\log\log z\ge\frac12\log z$. Theorem~1 is proved.

\medskip
The extension of Theorem~1 is actually a reformulation of Theorem~4.1 
in~[2], and its proof is worked out there in detail. Nevertheless, I 
briefly discuss this result to get a better understanding of it. 
Its proof is based on two propositions, and one of them is actually a 
weakened version of Theorem~1 of this paper.

\medskip\noindent
{\it On the proof of the extension of Theorem~1.}\/ This result 
is proved in~[2] with the help of two results formulated in 
Propositions~6.1 and~6.2 of that work. I discuss their content,
and show that Proposition~6.2 is a weakened version of 
Theorem~1. First I reformulate a slightly modified
version of it in the following Theorem~3.1.

\medskip\noindent
{\bf Theorem~3.1.} {\it Let us have a probability measure 
$\mu$ on a measurable space $(X,{\Cal X})$ together with a 
sequence of independent and $\mu$ distributed random 
variables $\xi_1,\dots,\xi_n$, $n\ge2$, and a countable, 
$L_1$-dense class $\Cal F$ of functions $f=f(x)$ on $(X,{\Cal X})$ 
with some parameter $D\ge1$ and exponent $L\ge1$ which 
satisfies the conditions $\sup\limits_{x\in X}|f(x)|\le 1$, 
$\int f(x)\mu(\,dx)=0$ and $\int f^2(x) \mu(\,dx)\le \sigma^2$ 
for all $f\in\Cal F$ with some $0<\sigma\le1$ such that the 
inequality $n\sigma^2>L\log n+\log D$ holds. Then there exists
a threshold index $A_0$ such that the normalized random sums
$S_n(f)$, $f\in {\Cal F}$, introduced in Theorem~1 satisfy the
inequality
$$
P\left(\sup_{f\in{\Cal F}}|S_n(f)|\ge A n^{1/2}\sigma^2\right)\le
e^{-A^{1/2}n\sigma^2/2}\quad \text{if } A\ge A_0. \tag3.8
$$
}

\medskip
I show that the estimate~(3.8) in Theorem~3.1 is a weakened
version of formula~(1.1) of Theorem~1. First I show that the
probability at the left-hand side of~(3.8) can be estimated 
by means of Theorem~1 in case~(c) with the choice 
$v=An^{1/2}\sigma^2$ if $A\ge A_0$ with a sufficiently large 
threshold index $A_0>0$.  We have to check that 
$v\ge u(\sigma)$ if $A_0$ is chosen sufficiently large. But
under the conditions of Theorem~3.1 
$n\sigma^2\ge L\log n\ge\frac18n\sigma^2$, and for 
$v=An^{1/2}\sigma^2$ we can write 
$v\ge\frac{A_0}{\sqrt n}n\sigma^2
\ge\frac{A_0}{2\sqrt n}n\sigma^2
+\frac{A_0}{2\sqrt n}(L\log n+\log D)
\ge\frac{C_5}{\sqrt n}(n\sigma^2+L\log n+\log D)=u(\sigma)$. 

Thus we can apply formula~(1.1) with $v=An^{1/2}\sigma^2$ 
to estimate the left-hand side of~(3.8), and we get the 
upper bound $C_1e^{-C_2\sqrt nv\log(v/\sqrt n\sigma^2)}=
C_1e^{-C_2 An\sigma^2\log A}$ if $A\ge A_0$. This is an
estimate sharper than formula~(3.8) if the (universal) 
constant~$A_0$ is chosen sufficiently large.  This 
calculation also indicates that Theorem~3.1 provides 
such a good estimate as Theorem~1 if $A\le\bar A_0$ with 
a fixed universal constant~$\bar A_0$. (Here we are not 
interested in the value of the universal constants in 
our estimates.) Besides, to prove the extension of 
Theorem~1 it is enough to have good estimates only 
in this case.

I shall only briefly discuss the content of Proposition~6.1 
in~[2], the other main ingredient in the proof of the 
extension of Theorem~1. Its proof is based on a classical 
method, called the chaining argument in the literature. 
It provides a sharp estimate for the tail distribution of 
the supremum of Gaussian random variables. But if we are 
interested in the tail distribution of the supremum of 
normalized partial sums of independent and identically 
distributed random variables, like in the extension of 
Theorem~1, then it only provides a weaker estimate. 
Proposition~6.1 actually contains the result we can get 
in our case with the help of the chaining method. This 
result is not sufficient for our purposes, but its 
combination with Theorem~3.1 enables us to prove the 
extension of Theorem~1. 

Here I do not discuss the details of the chaining argument.
It has a fairly detailed description in~[6], but also~[2] 
may help in understanding this method. I only remark that 
this method supplies a weaker estimate for the supremum 
of normalized partial sums of i.i.d. random variables, than 
for the supremum of the Gaussian random variables, because 
the tail distribution of partial sums of independent random 
variables has a slightly worse behaviour than the Gaussian 
tail distribution.

The main result of Proposition~6.1 in~[2] states that
under the conditions of the extension of Theorem~1 
there exists such a set of functions 
$\Cal F_{\bar\sigma}\subset\Cal F$ with some nice 
properties for which the inequality
$$
P\left(\sup_{f\in{\Cal F}_{\bar\sigma}} |S_n(f)|
\ge\frac u{\bar A}\right)\le 4\exp\left\{-\alpha\left(
\frac u{10\bar A\sigma}\right)^2\right\} \tag3.9
$$
holds if the number $u$ satisfies the condition
$n\sigma^2\ge(\frac u\sigma)^2\ge C_6(L\log\tfrac2\sigma+\log D)$ 
with a fixed constant $\bar A\ge1$, and the (sufficiently 
large) number $C_6=C_6(\bar A)$ appearing in the condition of 
formula~(3.9) depends on it. (The number~$\bar A$ was 
introduced in this estimate because of some technical 
reasons.) Moreover, the class of functions 
$\Cal F_{\bar\sigma}$, where $\bar\sigma=\bar\sigma(u)$ 
depends on the number $u$ in the above estimate has some 
properties which can be interpreted so that 
$\Cal F_{\bar\sigma}$ is a relatively small and dense subset 
of $\Cal F$. The set $\Cal F_{\bar\sigma}=\{f_1,\dots,f_m\}$ 
has $m\le D\bar\sigma^{-L}$ elements, and the sets 
$\Cal D_j=\{f\colon\;f\in\Cal F,\,\int (f-f_j)^2\,d\mu\le\bar\sigma^2\}$,
$1\le j\le m$, cover the set $\Cal F$, i.e. 
$\bigcup\limits_{j=1}^m\Cal D_j=\Cal F$. Theorem~6.1 also
provides some control on $\bar\sigma$. Namely, 
$\frac1{16}(\frac u{\bar A\bar\sigma})^2\ge n\bar\sigma^2
\ge\frac1{64}\left(\frac u{\bar A\sigma}\right)^2$, and 
the inequality $n\bar\sigma^2\ge L\log n+\log D$ also holds.
 
Formula~(3.9) gives a good (Gaussian type) estimate for a
supremum of partial sums. But in this estimate we took the
supremum for a class of functions 
$\Cal F_{\bar\sigma}\subset\Cal F$ instead of the class of 
functions $\Cal F$. The chaining argument does not enable 
us to give a good estimate if we take the supremum for a 
subclass of $\Cal F$ much larger than $\Cal F_{\bar\sigma}$. 
On the other hand, we get good estimates leading to the 
proof of the extension of Theorem~1 with the help of
Proposition~6.1 and a good bound on the probabilities 
$P\(\sup\limits_{h\in \Cal D'_j}S_n(h)\ge\frac u2\)$, 
$1\le j\le m$, where 
$\Cal D_j'=\{h=f-f_j\colon\; f_j\in\Cal D_j\}$.
Such a good bound can be obtained with the 
help of Theorem~3.1. But to apply this result we have to know 
that $n\bar\sigma^2\ge L\log n+\log D$, and this is the 
reason why this relation had to be proved in Proposition~6.1 
of~[2]. To get the desired estimate we also have to show that 
the number $m$ of the sets $\Cal D_j$ is not too large. Since
$m\le D\bar\sigma^{-L}$ this can be proved with the help of
the additional estimates on $\bar\sigma$ in Proposition~6.1. 
The proof of the extension of Theorem~1 with the help of 
Propositions~6.1 and~6.2 of~[2] is contained in that work, 
so I omit the details. They show some similarity to the
final step of the proof of Theorem~1 in this paper.

\medskip
I finish this paper with the formulation of some comments and 
problems.

\beginsection 4. Some comments on the methods and results of 
this paper.

\medskip\noindent
Our goal in this paper was to give a sharp estimate for the 
supremum of normalized partial sums $S_n(f)$, $f\in\Cal F$, 
of i.i.d. random variables in such cases that were not 
covered by previous results. The classical methods, like 
the chaining argument do not work in the study of such 
problems, since we have to bound some events in their 
applications for which we cannot give a sufficiently good
estimate. We have to deal with events whose probabilities 
are much larger than the value suggested by a Gaussian
comparison. I wanted to find a method that would work also 
in such cases.

One natural candidate for it was the so-called symmetrization 
argument. The extension of Theorem~1 presented in this paper 
was proved with the help of this method in~[2]. In the 
application of this method we consider a sequence of 
independent random variables $\e_1,\dots,\e_n$ with binomial 
distribution, i.e. $P(\e_j=1)=P(\e_j=-1)=\frac12$, 
$1\le j\le n$, which is independent also of the random 
variables $\xi_1,\dots,\xi_n$, and we reduce, with the help 
of some non-trivial inequalities,  the estimation of 
$P\left(\sup\limits_{f\in\Cal F}\frac1{\sqrt n}
\sum\limits_{j=1}^nf(\xi_j)>v\right)$
to the estimation of its `symmetrized version'
$P\left(\sup\limits_{f\in\Cal F}\frac1{\sqrt n}
\sum\limits_{j=1}^n\e_jf(\xi_j)>v\right)$.

My original plan was to prove Theorem~1 by means of a 
refinement of the symmetrization argument. But I met 
hard problems when I tried to carry out this program. 
The proof with the help of the symmetrization argument 
would have required the application of such an induction 
procedure, where at the start we need a good estimate 
for the supremum of the normalized partial sums $S_n(f)$ 
of very small random terms $f(\xi_j)$. More explicitly 
we should have handled the case when the expectation 
of the absolute value $E|f(\xi_j)|$ of the terms in the
sum are very small for all $f\in\Cal F$. In such cases the 
symmetrization argument is not useful, since if the terms
$f(\xi_j)$ in the normalized sum $S_n(f)$ are small, then 
the cancellation effect of the randomization by means of 
the multiplying factors $\e_j$, i.e. the replacement of 
the terms $f(\xi_j)$ by $\e_jf(\xi_j)$ is negligible. 
This implies that the symmetrization argument is 
ineffective in this case.

Hence a new method had to be found to estimate the supremum
of the normalized sums $S_n(f)$ if the additive terms
$f(\xi_j)$ in these sums are small. This was done in 
paper~[3]. After proving this result I recognized that
it makes the symmetrization argument in the study of 
the original problem superfluous. Theorem~1 can be proved 
in a much simpler direct way with the help of the result 
of~[3]. This was done in the present paper. 

The question arose for me whether the symmetrization
argument cannot be replaced by a simpler and stronger
method in the investigation of other problems. In 
particular, it would be interesting to consider the 
multivariate version of the extension of Theorem~1 
formulated in Theorem~8.4 of~[2]. This is an estimate
about the tail distribution of the supremum of 
appropriate degenerated $U$-statistics. This result 
was proved in~[2] by means of an adaptation of the 
symmetrization argument. We needed a multivariate 
version of this method which was based on a 
generalized version of the symmetrization lemma 
presented in Lemma~15.2 of~[2]. The proof of this 
lemma was not difficult, but in its application we 
have to estimate a rather complicated conditional 
probability (see formula~(15.3) in~[2]), and this 
made the proof of the above mentioned Theorem~8.4 
rather hard. It seems very probable that one can 
find a much simpler proof with the help of the 
method of the present paper.

\medskip\medskip\noindent
{\bf References.}

\medskip
\item{[1]} J. Koml\'os, P. Major, G. Tusn\'ady: An approximation
of partial sums of independent rv.'s and the sample DF. I.
Z. Wahrscheinlichkeitstheorie verw. Gebiete 32, 111--131 (1975)
\item{[2]} P. Major: On the estimation of multiple random integrals
and $U$-statistics {\it Lecture Notes in Mathematics} vol 2079
(Springer) Heidelberg New York Dordrecht London
\item{[3]} P. Major:  Sharp estimate on the supremum of a class
of partial sums of small i.i.d. random variables. submitted to
Electron. J. of Probab.
\item{[4]} H. P. McKean Jr.:{\it Stochastic integrals.}
Academic Press New York London (1969)
\item{[5]} M. Talagrand: New concentration inequalities in
product spaces. {\it Invent. Math.} {\bf 126}, 505--563 (1996) 
\item{[6]} M. Talagrand: {\it The general chaining.} Springer 
Monographs in Mathematics. \hfill\break 
Springer--Verlag, Berlin Heidelberg New York (2005)

\bye